\documentclass{optica-article}

\journal{opticajournal} 

\articletype{Research Article}

\usepackage{lineno}

\usepackage{amsopn}

\begin{document}

\title{{Distribution and Moments of a Normalized Dissimilarity Ratio for two Correlated Gamma Variables}}

\author{Elise Colin,\authormark{1} Razvigor Ossikovski,\authormark{2,*}}

\address{\authormark{1}DTIS-ONERA, University Paris-Saclay, 91123 Palaiseau, France\\
\authormark{2} LPICM, CNRS, Ecole Polytechnique, Institut Polytechnique de Paris, 91128 Palaiseau, France}

\email{\authormark{*}elise.colin@onera.fr}





\begin{abstract}
We consider two random variables $X$ and $Y$ following correlated Gamma distributions, characterized by identical scale and shape parameters and a linear correlation coefficient $\rho$. Our focus is on the parameter:
\[
  D(X,Y) = \frac{|X - Y|}{X + Y},
\]
which appears in applied contexts such as dynamic speckle imaging, where it is known as the \textit{Fujii index}. In this work, we derive a closed-form expression for the probability density function of $D(X,Y)$ as well as analytical formulas for its moments of order $k$.
Our derivation starts by representing $X$ and $Y$ as two correlated exponential random variables, obtained from the squared magnitudes of circular complex Gaussian variables. By considering the sum of $k$ independent exponential variables, we then derive the joint density of $(X,Y)$ when $X$ and $Y$ are two correlated Gamma variables. Through appropriate varable transformations, we obtain the theoretical distribution of $D(X,Y)$ and evaluate its moments analytically. These theoretical findings are validated through numerical simulations, with particular attention to two specific cases: zero correlation and unit shape parameter.
\end{abstract}



\section{Introduction}


Gamma, Rayleigh, and Nakagami distributions are fundamental in the statistical modeling of positive real-valued random variables \cite{nakagami1960m,papoulis1965random,simon2008digital}. They are extensively studied and applied in various domains, particularly in wireless telecommunications, where they characterize signal fading, propagation channel fluctuations, and antenna responses \cite{obeid2022distributions, reig2002bivariate}. Their ability to capture variations in received signal power due to multipath propagation and scattering effects makes them essential tools for analyzing and optimizing communication system performance.

Besides in telecommunications, these distributions also play a crucial role in speckle imaging, which is widely used in radar, laser, and ultrasound applications. In this context, they naturally emerge as probabilistic models for the amplitude and intensity of complex electromagnetic and acoustic waves. In particular, Gamma distributions are commonly employed to describe speckle field intensities in both electromagnetic and acoustic imaging, providing a robust statistical framework for analyzing scattered wave phenomena \cite{goodman1976some}.

In this paper, we consider two correlated Gamma-distributed random variables $X$ and $Y$, sharing the same shape and scale parameters, and linked by a correlation coefficient.

The comparison of two positive random variables $X$ and $Y$ can be carried out using various measures, such as the absolute value of the normalized difference 
\begin{equation}
  D(X,Y) = \frac{|X-Y|}{X+Y}.
\end{equation}

This quantity takes values in $[0,1]$ and exhibits scale invariance under $(X,Y) \mapsto (cX,cY)$ for any $c>0$. $D(X,Y)$ is widely used in dynamic speckle imaging, where its estimation is commonly referred to as the Fujii index \cite{fujii1985blood,fujii1987evaluation,pieczywek2017exponentially}. In this context, the above Normalized Dissimilarity Ratio plays a crucial role as a contrast measure in speckle correlation imaging, enabling the visualization of slow-moving flows by analyzing temporal fluctuations in scattered light. This approach is particularly effective for detecting very slow movements, where the speckle correlation remains high, corresponding to correlation coefficients close to 1. Despite its widespread use in experimental studies, no rigorous statistical framework has been established to analyze its distribution, nor to justify its applicability across different correlation regimes.

In this work, we derive the probability density function of $D(X,Y)$ in closed form and obtain explicit formulas for its moments of any integer order $m$. Besides these analytical results, a key contribution of our study is to provide a systematic reconstruction of the full sequence of transformations leading to the final expressions of the probability density function and moments of $D(X,Y)$. This includes, in particular, the derivation of the joint density of $(X,Y)$, two correlated Gamma-distributed variables, starting from the statistical correlation of two circular complex Gaussian variables, which serve as a fundamental representation of speckle fields. This step-by-step approach not only clarifies the statistical foundations linking speckle field properties to Gamma-distributed intensities but also explicitly establishes the connection between the correlation of a complex field and the behavior of the normalized contrast measure.

We begin by expressing the joint probability density function (PDF) of two correlated exponential random variables $(X,Y)$, which arise as the squared magnitudes of two circular complex Gaussian variables with identical variances. As an intermediate result, we retrieve the well-known relationship between the correlation coefficient of the Gamma-distributed intensities and the statistical correlation of the underlying complex speckle field (see\cite{song2016simulation}). 

In the next section, we extend this result to derive the joint PDF of two correlated Gamma-distributed variables, obtained as the sum of $k$ independent exponential variables. The parameter $k$ can be interpreted as the number of independent speckles integrated within a single measurement, a quantity of interest in various imaging applications.

Building upon these results, we first recover the density function of the ratio $X/Y$ before performing a final change of variables to obtain the analytical distribution of $D(X,Y)$ and its closed-form moment expressions for any order $m$.

In the final section, we validate our theoretical findings through numerical simulations, generating random samples of correlated Gamma-distributed variables by sampling from underlying Gaussian distributions. We then apply the derived transformations step by step to obtain realizations of $D(X,Y)$ in the general case. These numerical samples allow us to construct empirical histograms, which we systematically compare to the theoretical probability density function derived in closed form. Additionally, we compute empirical moments from the generated data and juxtapose them with the analytical expressions. Finally, we examine two notable special cases: $k=1$, corresponding to exponentially distributed variables, and the uncorrelated case $\rho=0$.

\section{Joint PDF of the squared magnitudes of two correlated complex circular normal variables}

We consider two circularly symmetric complex Gaussian random variables $Z_1$ and $Z_2$, each with zero mean and equal variance $\sigma_z^2$. Their joint distribution is fully characterized by their covariance matrix:

\begin{equation}
\Sigma = 
\begin{bmatrix}
\sigma_z^2 & \rho_z \sigma_z^2 \\
\rho_z \sigma_z^2 & \sigma_z^2
\end{bmatrix},
\end{equation}
where $\rho_z$ denotes the correlation coefficient between $Z_1$ and $Z_2$. Note that $\rho_z$  can be assumed to be real since $Z_1$ and $Z_2$ are circularly symmetric.

The inverse of $\Sigma$ being given by:

\begin{equation}
\Sigma^{-1} = \frac{1}{\sigma_z^2(1-\rho_z^2)}
\begin{bmatrix}
1 & -\rho_z \\
-\rho_z & 1
\end{bmatrix},    
\end{equation}
the PDF of $(Z_1, Z_2)$ is thus \cite{papoulis1965random}:

\begin{equation}
f_{Z_1,Z_2}(z_1, z_2) = \frac{1}{4\pi^2 \sigma_z^4 (1-\rho_z^2)} 
\exp\left( -\frac{1}{2\sigma_z^2(1-\rho_z^2)} 
\begin{bmatrix} z_1^* & z_2^* \end{bmatrix}
\begin{bmatrix}
1 & -\rho_z \\
-\rho_z & 1
\end{bmatrix}
\begin{bmatrix} z_1 \\ z_2 \end{bmatrix} 
\right).
\end{equation}
In polar coordinates $z_k = r_k e^{i\varphi_k}$, the PDF becomes:

\begin{equation}
f_{Z_1,Z_2}(r_1, r_2, \varphi_1, \varphi_2) =
\frac{r_1 r_2}{4\pi^2 \sigma_z^4 (1-\rho_z^2)}
\exp \left( -\frac{r_1^2 + r_2^2 - 2\rho_z r_1 r_2 \cos(\varphi_1 - \varphi_2)}{2\sigma_z^2 (1-\rho_z^2)} \right).
\end{equation}
By introducing the two random variables $X_1 = R_1^2$ and $X_2 = R_2^2$ for the squared magnitudes, we obtain:

\begin{equation}
f_{Z_1,Z_2}(x_1, x_2, \varphi_1, \varphi_2) = 
\frac{1}{16\pi^2 \sigma_z^4 (1-\rho_z^2)}
\exp \left( -\frac{x_1 + x_2 - 2\rho_z \sqrt{x_1 x_2} \cos(\varphi_1 - \varphi_2)}{2\sigma_z^2 (1-\rho_z^2)} \right).  
\end{equation}

To find the marginal distribution of $(X_1, X_2)$, we integrate over the angular variables:

\begin{equation}
f_{X_1,X_2}(x_1, x_2) = \int_0^{2\pi} \int_0^{2\pi} f_{Z_1,Z_2}(x_1, x_2, \varphi_1, \varphi_2) \, d\varphi_1 d\varphi_2.    
\end{equation}

To evaluate the double integral, we set $z = \frac{2\rho_z\sqrt{x_1 x_2}}{2\sigma_z^2(1-\rho_z^2)}$ in the integral definition of the modified Bessel function of zeroth order (3.339 in \cite{gradshteyn2014table}):

\begin{equation}
\int_0^{2\pi} \exp(z\cos \theta) d\theta = 2\pi I_0(z),
 \end{equation}
and obtain:

\begin{equation}
 \int_0^{2\pi}\int_0^{2\pi}\exp\left[z\cos(\varphi_1-\varphi_2)\right] d\varphi_1 d\varphi_2 = 4\pi^2 I_0\left(\frac{2\rho_z\sqrt{x_1 x_2}}{2\sigma_z^2(1-\rho_z^2)}\right).   
\end{equation}

Substituting this result in the double integral yields for the joint PDF of $(X_1, X_2)$:

\begin{equation}
f_{X_1,X_2}(x_1, x_2) = \frac{1}{(2\sigma_z^2)^2 (1-\rho_z^2)} \exp \left( -\frac{x_1 + x_2}{2\sigma_z^2 (1-\rho_z^2)} \right) I_0\left( \frac{2\rho_z \sqrt{x_1 x_2}}{2\sigma_z^2(1-\rho_z^2)} \right),
\end{equation}
confirming the well-known result that $(X_1, X_2)$ follow a correlated exponential distribution \cite{papoulis1965random}.

To derive the moments of the distribution, we make use of two key integrals. 
The starting point is the formula:

\begin{equation}
\int_0^{\infty} x^{\nu+1} e^{-\alpha x^2} J_{\nu}(\beta x) \, dx
= \frac{\beta^\nu}{(2\alpha)^{\nu+1}} \exp\left( -\frac{\beta^2}{4\alpha} \right),
\end{equation}
(6.631.4 in \cite{gradshteyn2014table}, valid for $\Re(\nu)>-1$ and $\Re(\alpha)>0$), in which we set $\nu=0$, substitute $i\beta\sqrt{x_2}$ for $\beta$, using the identity $J_{\nu}(ix) = i^{\nu} I_{\nu}(x)$ to switch from Bessel functions to modified Bessel functions, and perform the change of variable $x^2 = x_1$. This yields the first integral of interest:

\begin{equation}
\label{eq:firstIntegral}
\int_0^{\infty} e^{-\alpha x_1} I_0 \left( \beta \sqrt{x_1 x_2} \right) dx_1
= \frac{1}{\alpha} \exp\left( \frac{\beta^2}{4\alpha} x_2 \right).
\end{equation}

Taking the derivative with respect to $\alpha$ gives the second integral:

\begin{equation}
\label{eq:secondIntegral}
\int_0^{\infty} x_1 e^{-\alpha x_1} I_0 \left( \beta \sqrt{x_1 x_2} \right) dx_1
= \frac{1}{\alpha^2} \left( 1 + \frac{\beta^2}{4\alpha} x_2 \right) 
\exp\left( \frac{\beta^2}{4\alpha} x_2 \right).
\end{equation}

In the subsequent caclulations, we set $\alpha = \frac{1}{2\sigma_z^2(1-\rho_z^2)}$, $\beta = \frac{2\rho_z}{2\sigma_z^2(1-\rho_z^2)}$ and $C = \frac{1}{(2\sigma_z^2)^2(1-\rho_z^2)}$.

The first moment of $X_1$ (and, by symmetry, that of $X_2$ too) is found to be:

\begin{equation}
\langle X_1 \rangle = \langle X_2 \rangle = C\int_0^{2\pi}\int_0^{2\pi}x_2 e^{-\alpha (x_1+x_2)} I_{0}(\beta \sqrt{x_1 x_2}) \ dx_1 dx_2 = 2\sigma_z^2,
\end{equation}
by using the first integral of interest and an integration by parts.

Similarly, the second moment of $X_1$ is:

\begin{equation}
\langle X_1^2 \rangle = \langle X_2^2 \rangle = C\int_0^{2\pi}\int_0^{2\pi}x_2^2 e^{-\alpha (x_1+x_2)} I_{0}(\beta \sqrt{x_1 x_2}) \ dx_1 dx_2 =2(2\sigma_z^2)^2,
\end{equation}
after applying again the first integral and performing two integrations by parts. Consequently, the variance is:

\begin{equation}
\sigma^2_{X_1} = \sigma^2_{X_2} = \langle X_1^2 \rangle  - \langle X_1 \rangle^2 =(2\sigma_z^2)^2,
\end{equation}
which is consistent with the exponential distribution property that the mean equals the standard deviation.

The mixed second moment $\langle X_1 X_2 \rangle$ is obtained with the help of the second integral of interest and two integrations by parts:

\begin{equation}
 \langle X_1 X_2 \rangle = C\int_0^{2\pi}\int_0^{2\pi}x_1x_2 e^{-\alpha (x_1+x_2)} I_{0}(\beta \sqrt{x_1 x_2}) \ dx_1 dx_2 = (2\sigma_z^2)^2(1+\rho_z^2).   
\end{equation}

From the moments thus obtained, we find the correlation coefficient to be:

\begin{equation}
    \rho_{X_1, X_2} = \frac{\langle X_1 X_2 \rangle - \langle X_1 \rangle\langle X_2 \rangle}{\sigma_{X_1} \sigma_{X_2}} = \rho_z^2.
\end{equation}

\section{Derivation of the Joint PDF of Two Correlated Gamma Variables}

Having established the joint probability density function (PDF) of two correlated squared magnitudes following an exponential distribution, we now derive the joint PDF of two correlated Gamma-distributed variables. This follows naturally from the property that the sum of $ k $ independent exponentially distributed random variables follows a Gamma distribution with shape parameter $ k $.

We start from the previously obtained joint PDF of two correlated squared magnitudes, where we set $\sigma=\sigma_{X_1}=\sigma_{X_2}=2\sigma_z^2$ and $\rho=\rho_{X_1, X_2}=\rho_z^2$:

\begin{equation}
\label{eq:jointPdfSquaredMagnitudes}
f_{\text{exp}}(x_1, x_2) = \frac{1}{\sigma^2 (1 - \rho)}
\exp \left( -\frac{x_1 + x_2}{\sigma (1 - \rho)} \right)
I_0 \left( \frac{2\sqrt{\rho x_1 x_2}}{\sigma (1 - \rho)} \right),
\end{equation}

The above substitution emphasizes the fact that we are now dealing directly with the exponential distribution of the squared magnitudes $X_1$ and $X_2$, rather than with the complex fields $Z_1$ and $Z_2$ because $\sigma^2$ and $\rho$ are, respectively, the variance and the correlation coefficient of $X_1$ and $X_2$, as shown in the previous section. We now seek to extend this result to the sum of $ k $ independent correlated exponential variables, leading to correlated Gamma-distributed variables.

Following an approach analogous to that used in deriving the Nakagami distribution \cite{nakagami1960m}, we make use of the following two-dimensional Laplace transform formula:

\begin{align}
\frac{1}{\Gamma(k) b^{k-1}} 
\int_0^{\infty} \int_0^{\infty} (x_1 x_2)^{\frac{k-1}{2}} 
e^{-\alpha (x_1 + x_2)} I_{k-1} \left( 2b \sqrt{x_1 x_2} \right) 
e^{-z_1 x_1 - z_2 x_2} \,dx_1 dx_2 \notag \\
= \frac{1}{\left[ (z_1 + a)(z_2 + a) - b^2 \right]^k}.
\end{align}

This formula, reported as equation (79) in \cite{voelker2013zweidimensionale}, holds for $\Re(k) > 0$. By setting the parameters to match our exponential distribution, $
k=1$, $\alpha = \frac{1}{\sigma(1-\rho)}$, $b = \frac{\sqrt{\rho}}{\sigma(1-\rho)}
$, 
we obtain the joint characteristic function of the exponential distribution:

\begin{equation}
\label{eq:jointCharFunc}
\varphi_{\text{exp}}(z_1, z_2) =
\frac{1}{\sigma^2 (1 - \rho) \left[ (z_1 + a)(z_2 + a) - b^2 \right]}.
\end{equation}

Since the sum of $ k $ independent exponential variables follows a Gamma distribution, the characteristic function of two correlated Gamma variables can be directly deduced:
\begin{equation}
\label{eq:charFuncCorrelatedGamma}
\varphi_{\text{Gamma}}(z_1, z_2) = \left[ \varphi_{\text{exp}}(z_1, z_2) \right]^k 
= \frac{1}{\sigma^{2k} (1 - \rho)^k \left[ (z_1 + a)(z_2 + a) - b^2 \right]^k}.
\end{equation}

To invert this characteristic function and recover the probability density function, we use the two-dimensional Mellin inversion formula \cite{Nicolas2019Mellin}:

\begin{align}
\left( \frac{1}{2\pi i} \right)^2 
\int_{c - i\infty}^{c + i\infty} \int_{c - i\infty}^{c + i\infty} 
\frac{1}{\left[ (z_1 + a)(z_2 + a) - b^2 \right]^k} 
e^{z_1 x_1 + z_2 x_2} \,dz_1 dz_2 \notag \\
= \frac{1}{\Gamma(k) b^{k-1}} (x_1 x_2)^{\frac{k-1}{2}} 
e^{-\alpha (x_1 + x_2)} I_{k-1} \left( 2b \sqrt{x_1 x_2} \right).
\end{align}
Applying this inversion formula to the characteristic function in \eqref{eq:charFuncCorrelatedGamma}, we finally obtain the joint PDF of two correlated Gamma-distributed variables:

\begin{equation}
\label{eq:jointPdfCorrelatedGamma}
f_{\text{Gamma}}(x_1, x_2) =
\frac{(x_1 x_2)^{\frac{k-1}{2}}}
{\Gamma(k) \sigma^{k+1} (1 - \rho) \rho^{\frac{k-1}{2}}}
\exp \left( -\frac{x_1 + x_2}{\sigma (1 - \rho)} \right)
I_{k-1} \left( \frac{2\sqrt{\rho x_1 x_2}}{\sigma (1 - \rho)} \right).
\end{equation}
This expression generalizes the previous result for the exponential PDF with $k=1$ and agrees with results reported in the literature (see \cite{bithas2007distributions, piboongungon2005bivariate}).

\section{Probability Density Function of the Ratio of two correlated Gamma-distributed random variables}

In this section, we derive the PDF of the ratio $ Z = X_1 / X_2 $, where $ X_1 $ and $ X_2 $ are two correlated Gamma-distributed random variables with the same shape and scale parameters.

By definition, the PDF of the ratio $ Z $ for two non-negative random variables $ X_1 $ and $ X_2 $ is given by:

\begin{align}
\label{eq:pdfRatioDef}
f_Z(z) &= \int_0^{\infty} x f(zx, x) \,dx.
\end{align}
By substituting the joint PDF $ f_{\text{Gamma}}(x_1, x_2) $ derived previously, we have:

\begin{align}
f_Z(z) &= \int_0^{\infty} x \frac{(zx^2)^{\frac{k-1}{2}}}
{\Gamma(k) \sigma^{k+1} (1 - \rho) \rho^{\frac{k-1}{2}}}
\exp \left( -\frac{x (z+1)}{\sigma (1 - \rho)} \right)
I_{k-1} \left( \frac{2x \sqrt{\rho z}}{\sigma (1 - \rho)} \right) dx.
\end{align}
Factoring out terms independent of $ x $ and simplifying yields:

\begin{align}
f_Z(z) &= \frac{z^{\frac{k-1}{2}}}{\Gamma(k) \sigma^{k+1} (1 - \rho) \rho^{\frac{k-1}{2}}}
\int_0^{\infty} x^k 
\exp \left( -\frac{z+1}{\sigma (1 - \rho)} x \right)
I_{k-1} \left( \frac{2\sqrt{\rho z}}{\sigma (1 - \rho)} x \right) dx.
\end{align}
To solve the integral, we use the formula (6.623.2 in \cite{gradshteyn2014table}):

\begin{equation}
\label{eq:formula66232}
\int_{0}^{\infty} x^{\nu+1} e^{-\alpha x} J_{\nu} (\beta x) \,dx 
= \frac{2\alpha (2\beta)^{\nu} \Gamma \left( \nu + \frac{3}{2} \right)}
{\sqrt{\pi} \left( \alpha^2 + \beta^2 \right)^{\nu + \frac{3}{2}}},
\end{equation}
valid for $\Re(\nu)>-1 $ and $ \Re(\alpha)>|\Im(\beta)|$. In our case, we set:
$\nu = k-1$,  $\alpha = \frac{z+1}{\sigma(1-\rho)}$, and $ \beta = i\frac{2\sqrt{\rho z}}{\sigma(1-\rho)}$,
and make use of the identity $ J_{\nu}(ix)=i^{\nu}I_{\nu}(x) $. Application of this formula yields:

\begin{equation}
f_Z(z) = \frac{\Gamma\left(k + \frac{1}{2}\right) z^{\frac{k-1}{2}} 
\left( \frac{z+1}{\sigma(1-\rho)} \right) 
\left( \frac{2\sqrt{\rho z}}{\sigma(1-\rho)} \right)^{k-1}}
{\Gamma(k) \sqrt{\pi} \sigma^{k+1} (1-\rho) \rho^{\frac{k-1}{2}} 
\left[ \left( \frac{z+1}{\sigma(1-\rho)} \right)^2 - \left( \frac{2\sqrt{\rho z}}{\sigma(1-\rho)} \right)^2 \right]^{k+\frac{1}{2}}}.
\end{equation}
To simplify further, we use Legendre’s doubling formula for the Gamma function (8.335.1 in \cite{gradshteyn2014table}):

\begin{equation}
\Gamma \left( k + \frac{1}{2} \right) =
\frac{\sqrt{\pi} \, \Gamma(2k)}{2^{2k-1} \Gamma(k)},
\end{equation}
and the well-known relationship between the Beta and Gamma functions (8.384.1 in \cite{gradshteyn2014table}):

\begin{equation}
B(\alpha, \beta) = \frac{\Gamma(\alpha) \Gamma(\beta)}{\Gamma(\alpha + \beta)}.
\end{equation}
These substitutions lead to the final expression for the PDF of $ Z $:

\begin{equation}
\label{eq:pdfZfinal}
f_Z(z) = \frac{(1-\rho)^k z^{k-1} (z+1)}{B(k,k) \left[ (z+1)^2 - 4\rho z \right]^{k+\frac{1}{2}}}.
\end{equation}

This result is a generalization of the known expression for two correlated exponential variables \cite{malik1986probability} and is a special case of that for two correlated generalized Gamma variables \cite{bithas2007distributions}.

\section{Probability Density Function of the Normalized Dissimilarity Ratio}

The goal of this derivation is to determine the PDF of the Normalized Dissimilarity Ratio:
\begin{equation}
D = \frac{|X_1 - X_2|}{X_1 + X_2}.
\end{equation}
As already mentioned, this parameter is widely used in dynamic speckle imaging, where it is commonly referred to as the Fujii index. To obtain its distribution, we apply a change of variable from the ratio $Z = X_1 / X_2$ to $D$.

We introduce the auxiliary variable:
\begin{equation}
D' = \frac{Z-1}{Z+1},
\end{equation}
which provides a one-to-one mapping between $Z$ and $D'$, with the inverse relation:
\begin{equation}
Z = \frac{1 + D'}{1 - D'}.
\end{equation}
 The standard formula for a change of variable gives:

\begin{align}
\label{eq:pdfRprime}
f_{D'}(r') &= f_Z \left( \frac{1 + r'}{1 - r'} \right) 
\left| \frac{d}{dr'} \left( \frac{1 + r'}{1 - r'} \right) \right|.
\end{align}
Substituting the previously derived PDF of $Z$, we get:

\begin{align}
f_{D'}(r') &= \frac{(1 - \rho)^k \left( \frac{1 + r'}{1 - r'} \right)^{k-1} }
{B(k,k) \left[ \left( \frac{2}{1 - r'} \right)^2 - 4\rho \left( \frac{1 + r'}{1 - r'} \right) \right]^{k+\frac{1}{2}}} 
\left| \frac{2}{(1 - r')^2} \right|. 
\end{align}
After simplification, we obtain:

\begin{equation}
f_{D'}(r') = \frac{(1 - \rho)^k (1 - r'^2)^{k-1}}{B(k,k) 2^{2k-1} (1 - \rho + \rho r'^2)^{k+\frac{1}{2}}}.
\end{equation}

Since $D = |D'|$, its PDF is given by summing the contributions from both positive and negative values of $D'$:

\begin{equation}
f_D(r) = f_{D'}(r') + f_{D'}(-r') = 2 f_{D'}(r'),
\end{equation}
which leads to the final expression:

\begin{equation}
f_D(r) = \frac{(1 - \rho)^k (1 - r^2)^{k-1}}
{B(k,k) 2^{2k-2} (1 - \rho + \rho r^2)^{k+\frac{1}{2}}}.
\end{equation}
for the PDF of the Normalized Dissimilarity Ratio.

\section{Analytical Expressions of the \texorpdfstring{$m$}{m}-th Moment of the Normalized Dissimilarity Ratio Distribution}

In this section, we derive closed-form expressions for the $m$-th order moments of the Normalized Dissimilarity Ratio $D$. The derivation uses integral transformations and properties of the hypergeometric function to yield three equivalent analytical formulations, each offering different computational advantages. 

\subsection*{First Formulation}

By definition, the $m$-th moment of $D$ is given by:
\begin{align}
\label{eq:mthMomentFujiiDef}
\langle D^m \rangle &= \int_0^1 r^m f_D(r) \, dr \notag \\
&= \int_0^1 r^m \frac{(1 - \rho)^k (1 - r^2)^{k-1}}{B(k,k) 2^{2k-2} (1 - \rho + \rho r^2)^{k+\frac{1}{2}}} \, dr.
\end{align}

By introducing the substitutions $t = r^2$ and $z = -\frac{\rho}{1 - \rho}$, and using the integral representation of the hypergeometric function (9.111 in \cite{gradshteyn2014table}):
\begin{equation}
\label{eq:hypergeomRepresentation}
F(\alpha, \beta; \gamma, z) = \frac{1}{B(\beta, \gamma - \beta)} \int_0^1 t^{\beta - 1} (1 - t)^{\gamma - \beta - 1} (1 - tz)^{-\alpha} \, dt,
\end{equation}
we obtain the first closed-form expression for the $m$-th moment:
\begin{equation}
\label{eq:hypergeomApplied}
\langle D^m \rangle = \frac{B \left( \frac{m+1}{2}, k \right) (1 - \rho)^{-\frac{1}{2}}}{B(k,k) 2^{2k-1}} F \left( k + \frac{1}{2}, \frac{m+1}{2}; k + \frac{m+1}{2}; -\frac{\rho}{1 - \rho} \right).
\end{equation}

\subsection*{Second Formulation}

Using the hypergeometric function transformation (9.131.1 GA 218 (92) in \cite{gradshteyn2014table}):
\begin{equation}
\label{eq:transformationFormula1}
F(\alpha, \beta; \gamma; z) = (1 - z)^{-\beta} F \left( \beta, \gamma - \alpha; \gamma; \frac{z}{z - 1} \right),
\end{equation}
and its symmetry property $F(\alpha, \beta; \gamma; z) = F(\beta, \alpha; \gamma; z)$, we arrive at an alternative formulation:
\begin{equation}
\label{eq:mthMomentFujiiFinal1}
\langle D^m \rangle = \frac{B \left( \frac{m+1}{2}, k \right) (1 - \rho)^{\frac{m}{2}}}{B(k,k) 2^{2k-1}} F \left( \frac{m}{2}, \frac{m+1}{2}; k + \frac{m+1}{2}; \rho \right).
\end{equation}

\subsection*{Third Formulation}

Finally, by employing another transformation (9.131.1 GA 218 (91) in \cite{gradshteyn2014table}):
\begin{equation}
\label{eq:transformationFormula2}
F(\alpha, \beta; \gamma; z) = (1 - z)^{-\alpha} F \left( \alpha, \gamma - \beta; \gamma; \frac{z}{z - 1} \right),
\end{equation}
we derive a third expression for the $m$-th moment:
\begin{equation}
\label{eq:mthMomentFujiiFinal2}
\langle D^m \rangle = \frac{B \left( \frac{m+1}{2}, k \right) (1 - \rho)^k}{B(k,k) 2^{2k-1}} F \left( k, k + \frac{1}{2}; k + \frac{m+1}{2}; \rho \right).
\end{equation}

These three formulations are mathematically equivalent and provide flexibility depending on the numerical context  \cite{simon2008digital}. The hypergeometric function ${}_2F_1$ is implemented in popular computational tools such as \texttt{MATLAB} and \texttt{Python}'s \texttt{scipy.special} library, enabling straightforward numerical evaluation of these moments. This versatility is particularly useful when analyzing the statistical properties of the Fujii index under various parameter settings.

\section{Numerical Validation and Discussion of the Derived Expressions}

In this section, we validate the closed-form expressions derived for the probability density function and moments of the Normalized Dissimilarity Ratio $D(X,Y) = \frac{|X - Y|}{X + Y} $. To ensure reproducibility and facilitate the exploration of different parameter settings, all computations were implemented in a shared Python notebook (Google Colab). We generate numerical samples of the random variables correlated with Gamma distributions $X$ and $Y$ first by sampling from correlated circular complex Gaussian distributions with predefined parameters $\sigma$ and $\rho$, and then applying the transformations leading to Gamma-distributed intensities with $k$ shape parameter.

The empirical distributions obtained from these generated samples allow us to construct histograms that are then compared to the theoretical probability density functions derived in the previous sections. Furthermore, empirical moments are computed from the generated data and juxtaposed with the closed-form expressions established analytically. This numerical approach not only confirms our theoretical results but also provides insights into the behavior of the Normalized Dissimilarity Ratio under various parameter configurations.

\subsection{Empirical Validation of the Correlation Structure}

In this section, we aim to verify numerically a key intermediate result that forms the basis of our approach: two exponentially distributed variables $X$ and $Y$ with correlation $\rho^2$ can be simply obtained from two correlated circular complex Gaussian variables $Z_x$ and $Z_y$.
, using the relation:
\[
Z_y = \rho Z_x + \sqrt{1 - \rho^2} W,
\]
where $W$ is an independent complex Gaussian variable. This ensures that the resulting intensities $X = |Z_x|^2$ and $Y = |Z_y|^2$ are correlated with a coefficient $\rho^2$. This verification step is crucial, as the entire derivation of the Normalized Dissimilarity Ratio distribution relies on this initial construction. It is worth noting that this result has already been established in \cite{song2016simulation}

\begin{figure}[ht]
    \centering
    \includegraphics[width=0.8\textwidth]{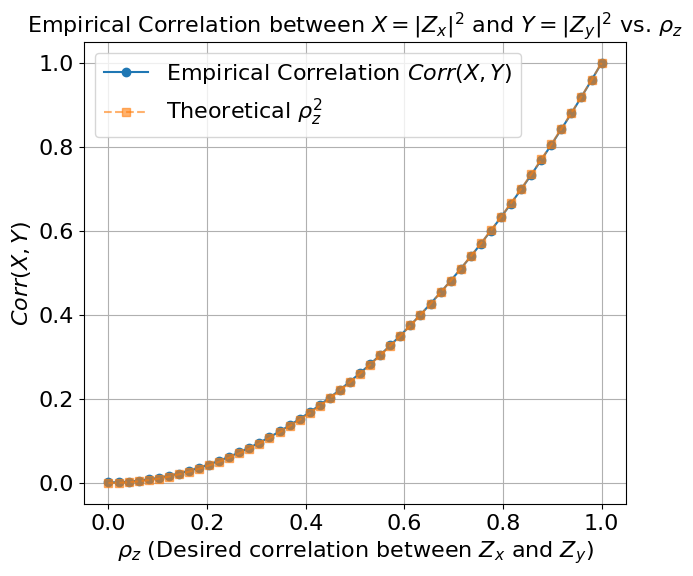}
    \caption{Empirical correlation between the Gamma-distributed intensities $ X = |Z_x|^2 $ and $ Y = |Z_y|^2 $ as a function of the input correlation $ \rho $ of the complex fields. The empirical values superimpose on the theoretical curve $ \rho^2 $, confirming the expected correlation structure.}
    \label{fig:correlation_relation}
\end{figure}

In the first notebook section, we reproduce this result empirically. We generate $N$ samples of $Z_x$ and $Z_y$ for different values of $\rho$, and compute the empirical correlation between $X$ and $Y$. Figure~\ref{fig:correlation_relation} shows that the empirical correlation closely matches the theoretical value $\rho^2$.

\subsection{Joint Probability Density of \texorpdfstring{$ (X,Y) $}{(X,Y)}}

In this section, we validate the derived theoretical joint probability density functions of $ (X,Y) $ by comparing them with empirical distributions obtained through numerical simulations. The variables $ X $ and $ Y $ are generated from correlated circular complex Gaussian fields, and their joint distributions are analyzed under different parameter settings.
We provide three types of visualizations to illustrate our theoretical results:
\begin{itemize}
    \item \textbf{Exponential Case}: Comparison between the theoretical and empirical joint density of $ (X,Y) $ when $ X $ and $ Y $ follow correlated exponential distributions (i.e., $ k=1 $), with a given set of parameters $ \rho $ and $ \sigma $.
    \item \textbf{Gamma Case}: Similar comparison for correlated Gamma-distributed variables $ (X,Y) $ with a fixed $ \rho $, $ \sigma $, and shape parameter $ k $.
    \item \textbf{3D Visualization}: A three-dimensional surface plot of the theoretical joint densities for both exponential and Gamma cases, highlighting the joint behavior and the associated marginal distributions.
\end{itemize}

\begin{figure}[ht]
    \centering
    \includegraphics[width=0.8\textwidth]{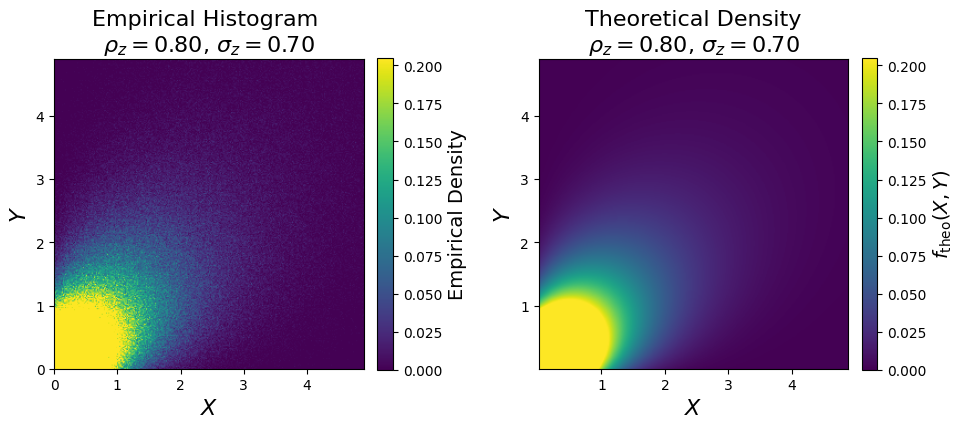}
    \caption{Comparison between empirical and theoretical joint probability density of $ (X,Y) $ in the exponential case ($ k=1 $) for $ \rho_z=0.8 $ and $ \sigma_z=0.7 $.}
    \label{fig:joint_pdf_exponential}
\end{figure}

For the empirical histograms presented, the example parameter values used are $ \rho_Z = 0.8 $ (corresponding to $ \rho = 0.64 $), $ \sigma_Z = 0.7 $ (hence $ \sigma = 2.88 $), and $ k = 12 $ in the case of Gamma-distributed variables. Figure~\ref{fig:joint_pdf_exponential} displays the empirical 2D histogram and the theoretical joint density contours for $ (X,Y) $ in the exponential case. The near-perfect overlap between the empirical histogram and the theoretical contours shows the consistency of our derived expression for the exponential joint PDF.

\begin{figure}[ht]
    \centering
    \includegraphics[width=0.8\textwidth]{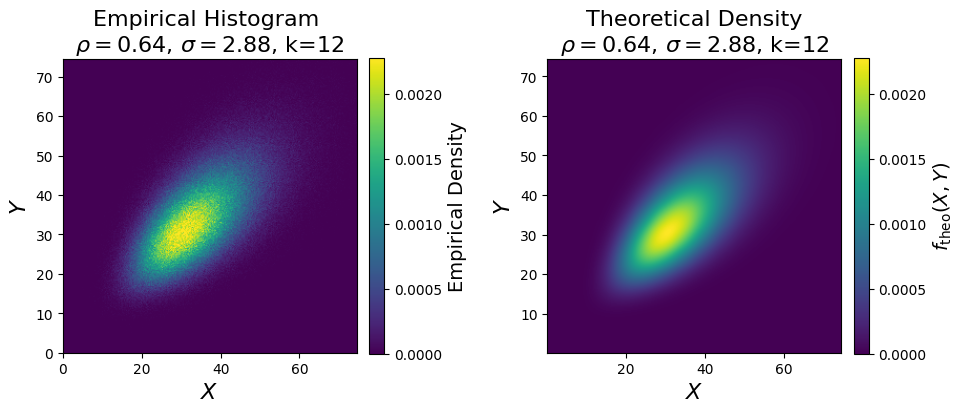}
    \caption{Comparison between empirical and theoretical joint probability density of $ (X,Y) $ in the Gamma case for $ \rho_z=0.8, \rho=0.64 $, $ \sigma_z=0.7,\sigma=2.88 $ and $k=12$}
    \label{fig:joint_pdf_gamma}
\end{figure}

In Figure~\ref{fig:joint_pdf_gamma}, we present a similar comparison for the Gamma-distributed case, with $ k > 1 $. Again, the agreement between empirical data and the theoretical model confirms the accuracy of our derivations, even for higher shape parameters and more complex correlation structures.

\begin{figure}[htbp]
    \centering
    \includegraphics[width=0.48\textwidth]{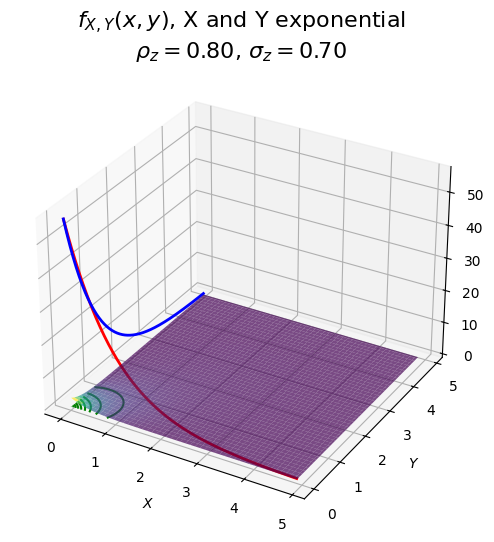}
    \includegraphics[width=0.48\textwidth]{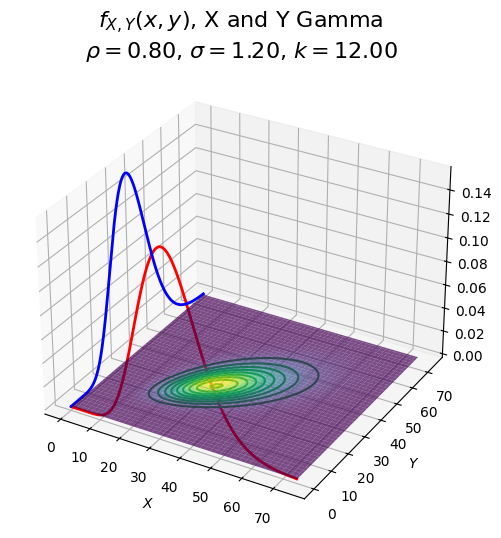}
    \caption{Visualization of the theoretical joint probability density of $ (X,Y) $ for the Exponential case (left) and the Gamma case (right). The joint PDFs are represented in a horizontal plane, while the third dimension is used to display the corresponding marginal distributions along the $X$- and $Y$-axes.}
    \label{fig:joint_pdf_3d}
\end{figure}

Finally, to better illustrate the relationship between the joint distribution and its marginals for both the exponential and Gamma cases, Figure~\ref{fig:joint_pdf_3d} presents a graphical representation in which the theoretical joint PDFs are displayed in a horizontal plane, while the third dimension is used specifically to depict the corresponding marginal distributions along the $X$- and $Y$-axes.

The corresponding section of the notebook ensures the agreement between theoretical derivations and empirical results by systematically varying the parameters $\sigma$, $\rho$, and $k$. This step-by-step approach enables us to check the accuracy of the joint probability density functions for both correlated exponential and Gamma-distributed variables. It also provides an insightful visualization of how these parameters influence the resulting distributions.

\subsection{Empirical Distribution of the Normalized Dissimilarity Ratio}

\begin{figure}[htbp]
    \centering
    \includegraphics[width=0.45\textwidth]{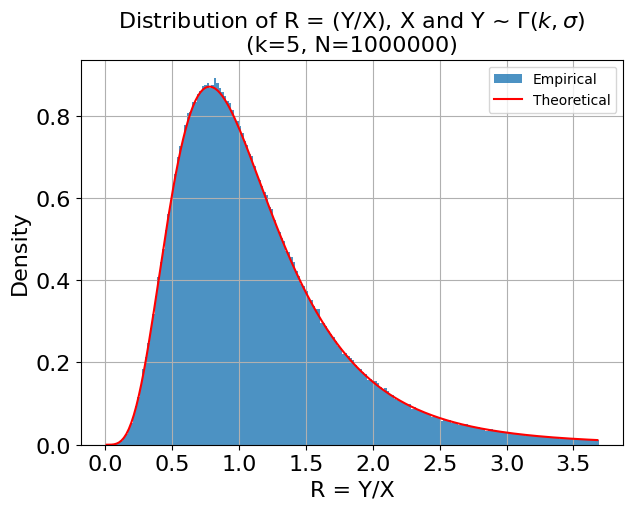}
    \includegraphics[width=0.45\textwidth]{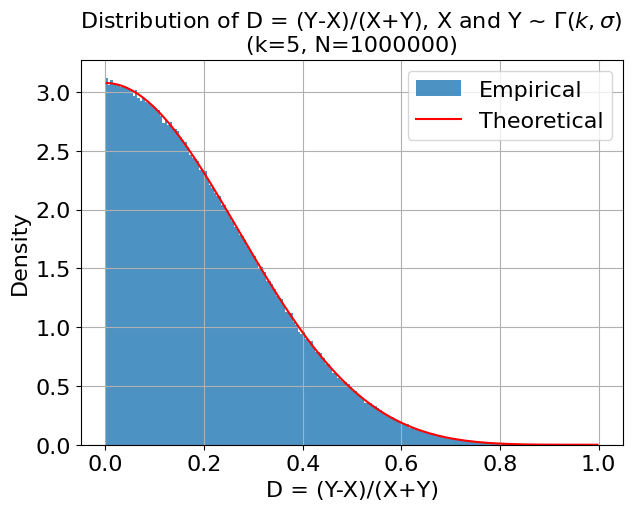}
    \caption{Comparison between the empirical histogram and the theoretical density of the Normalized Dissimilarity Ratio $ D(X,Y) $ for fixed values of $ \rho $, $ \sigma $, and $ k $.}
    \label{fig:fujii_pdf}
\end{figure}

In this subsection, we present the empirical histograms of the ratio $X/Y$ and the normalized dissimilarity parameter $D(X,Y)$, along with their corresponding theoretical probability density functions, in Figure \ref{fig:fujii_pdf}. While the results shown in the article are based on a single set of parameters of $\rho$, $\sigma$, and $k$, the accompanying notebook allows for comprehensive exploration by varying these parameters. Thus, the agreement between empirical simulations and theoretical expressions can be checked across a wide range of parameter values.

\subsection{Evolution of the Normalized Dissimilarity Ratio with \texorpdfstring{$ \rho $}{rho}}

\begin{figure}[htbp]
    \centering
    \includegraphics[width=0.45\textwidth]{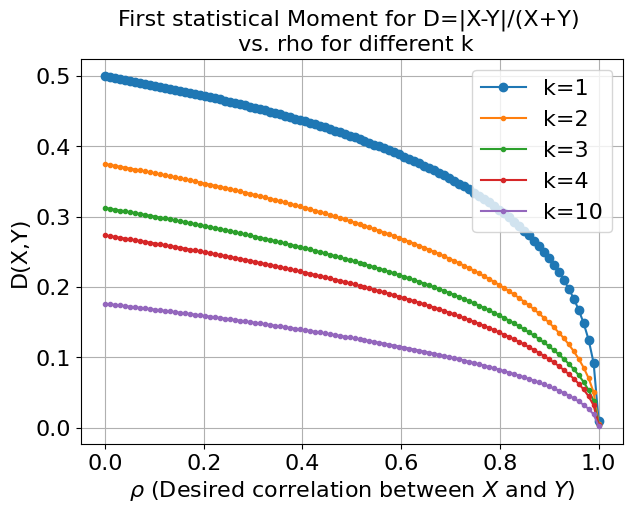}
    \includegraphics[width=0.45\textwidth]{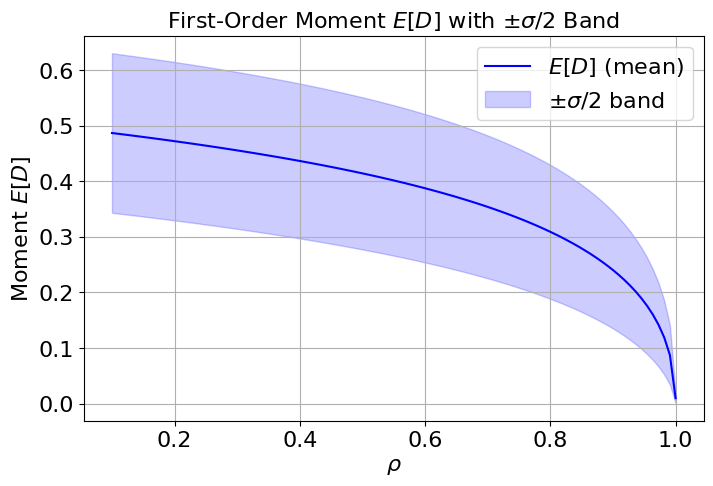}
    \caption{Evolution of the Normalized Dissimilarity Ratio index as a function of $\rho$ for different values of $k$.}
    \label{fig:fujii_vs_rho}
\end{figure}

Figure~\ref{fig:fujii_vs_rho} presents the evolution of the Normalized Dissimilarity Ratio as a function of the correlation coefficient $\rho$ for different values of $k$. As expected, the Normalized Dissimilarity Ratio decreases towards zero as $\rho$ increases, particularly for high correlation values. The most informative region of the curve is for values of $\rho$ close to 1, where the Normalized Dissimilarity Ratio exhibits significant sensitivity to small changes in correlation.

In the particular case of $k=1$ (exponential distributions), the initial value of $D$ at $\rho=0$ is 1/2, with $D$ \textbf{ being uniformly distributed between 0 and 1.} Additionally, we have included a plot representing the evolution of the mean value of $D$ as well as a band corresponding to $\pm 1/2$ of the standard deviation.

\subsection{Evolution of the Normalized Dissimilarity Ratio with \texorpdfstring{$ k $}{k} in the Uncorrelated Case}

\begin{figure}[htbp]
    \centering
    \includegraphics[width=0.8\textwidth]{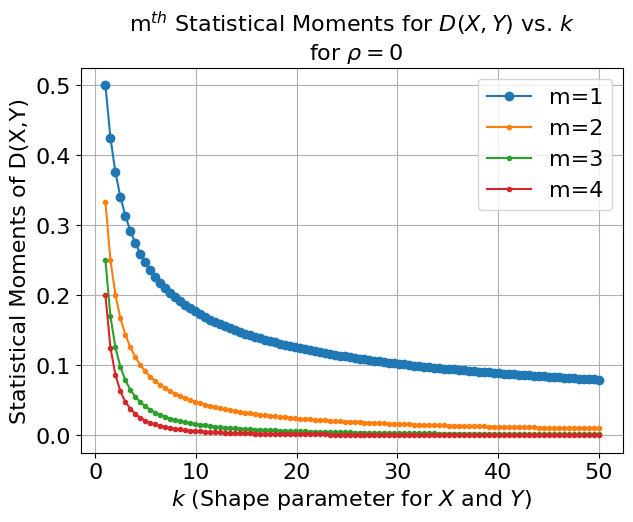}
    \caption{Evolution of the first four moments of the Normalized Dissimilarity Ratio as a function of $ k $ for $ \rho = 0 $.}
    \label{fig:fujii_vs_k}
\end{figure}

Finally, in the specific case of uncorrelated intensities $( \rho = 0)$, we investigate the dependence of the first four moments of Normalized Dissimilarity Ratio on the shape parameter $k$. Figure~\ref{fig:fujii_vs_k} shows that the Normalized Dissimilarity Ratio decreases linearly with $k$. This behavior confirms that increasing $k$ corresponds to increasing the number of independent speckles integrated into an image, reducing contrast variations reflected in the Normalized Dissimilarity Ratio.

\section{Conclusions}
\label{sec:conclusions}

In this work, we have derived closed-form expressions for the moments of the Normalized Dissimilarity Ratio $D(X,Y)=\left| X - Y \right| /(X+Y)$ when $X$ and $Y$ follow correlated Gamma distributions. Our approach stems from the representation of $X$ and $Y$ as squared magnitudes of correlated circular complex Gaussian variables, leading to the joint probability density function of $(X,Y)$ as correlated Gamma-distributed intensities. The expressions derived depend on four key parameters: the shape parameter $k$ of the Gamma distribution, the variance $\sigma$ of the intensity variables, the correlation coefficient $\rho$ (equal to the squared correlation coefficient of the underlying complex fields) and the moment order $m$.

We have explored two notable special cases. In the case of zero correlation $\rho=0$, the intensities correspond to fully decorrelated speckles. In this regime, the Normalized Dissimilarity Ratio decreases linearly as the shape parameter $k$ increases. This result implies a relationship between the Fujii index and speckle activity: for a fixed integration time, higher medium velocities increase the number of speckles integrated within an image, thereby reducing the observed normalized contrast. Similarly, increasing the number of grains per speckle reduces the Fujii index, highlighting its dependence on speckle properties and motion dynamics.

In the case of non-zero correlation with $k=1$, we recover the scenario of a fully developed speckle, such as a static speckle field with one grain per pixel. In this regime, the Fujii index ranges from 0 (for zero velocity and maximum correlation $\rho=1$) to $1/2$ as the correlation tends to zero. This case is particularly relevant for extremely slow movements, where the speckle correlation remains close to 1, making the Fujii index a sensitive tool for detecting very small motions.

Our results provide a deeper understanding of the parameterization required for dynamic speckle measurements, particularly in choosing the appropriate integration times, speckle grain sizes, and expected correlation levels. Beyond dynamic speckle imaging, the Normalized Dissimilarity Ratio and its moments, as derived in this work, could be applied to a wide range of contexts, including multichannel polarimetric speckle correlations and contrast analysis between correlated Gamma-distributed variables in broader statistical frameworks.

\section*{Supplementary Material}
The Python codes used to generate the figures in this article, as well as to verify our analytical formulas, are available as a Jupyter Notebook on Google Colab. This interactive environment allows for the direct execution and reproduction of our computational results. The notebook can be accessed at the following link:
\href{https://colab.research.google.com/drive/1wvGKSS6JPpl8BSzy9uykqFnWYLyzZjgJ\string?_usp=sharing}{Google Colab}

\section*{Acknowledgments}
We would like to thank our esteemed colleagues Enrique Garcia Caurel, Aurélien Plyer and Xavier Orlik for their invaluable participation in all our discussions and for their constant support throughout this work.

\bibliography{sample}
\end{document}